\documentstyle[amssymb]{article}

\newcommand{\proof}{{\bf Proof:  }}
\newcommand{\remark}{{\bf Remark:  }}
\newcommand{\remarks}{{\bf Remarks:  }}
\newcommand{\example}{{\bf Example:  }}
\newcommand{\examples}{{\bf Examples:  }}
\newcommand{\dimv}{\underline{\dim}\,}
\newcommand{\hb}{\newline\hspace*{\fill}$\Box$}

\newcommand{\ses}[3]
{\mbox{$0 \rightarrow #1 \rightarrow #2 \rightarrow #3 \rightarrow 0$}}
\newtheorem{theorem}{Theorem}[section]
\newtheorem{lemma}[theorem]{Lemma}
\newtheorem{definition}[theorem]{Definition}
\newtheorem{proposition}[theorem]{Proposition}
\newtheorem{corollary}[theorem]{Corollary}

\begin{document}
\title{\bf Quivers, desingularizations and canonical bases}
\author{Markus Reineke\\[2ex] BUGH Wuppertal, Gau\ss str. 20, D-42097
Wuppertal, Germany\\ (e-mail: reineke@math.uni-wuppertal.de)}
\date{}

\maketitle

\begin{abstract} A class of desingularizations for orbit closures of 
representations of Dynkin quivers is constructed, which can be viewed as a 
graded analogue of the Springer resolution. A stratification of the singular 
fibres is introduced; its geometry and combinatorics are studied. Via the Hall 
algebra approach, these constructions relate to bases of quantized enveloping 
algebras. Using Ginzburg's theory of convolution algebras, the base change 
coefficients of Lusztig's canonical basis are expressed as decomposition 
numbers of certain convolution algebras.
\end{abstract}

\parindent0pt

\section{Introduction}

The varieties of representations of Dynkin quivers are of central importance 
to geometric realizations of quantized enveloping algebras and structures 
related to them. For example, positive parts of quantized enveloping algebras 
can be realized by a convolution construction (C. M. Ringel's Hall algebra 
approach \cite{RiH}). Of particular interest are orbit closures in these 
varieties, since their intersection cohomology realizes G. Lusztig's canonical 
bases (\cite{L1}). On the other hand, these varieties can be viewed as 
(quiver-) graded analogues of the nilpotent cones, which makes their geometric 
analysis interesting in itself.\\[1ex]
Motivated by this last analogy, this paper starts a program to develop an 
analogue of Springer theory of nilpotent cones in the quiver context, and to 
explore the quantum group theoretic consequences of such constructions.\\[1ex]
As a first step, a model for Springer's resolution has to be found. Therefore, 
we construct in section 2 desingularizations of orbit closures for 
representations of arbitrary Dynkin quivers, generalizing the case of 
equioriented quivers of type $A$ treated in \cite{ADK}. At the heart of this 
construction lies the Auslander-Reiten theory of finite dimensional 
algebras.\\[1ex]
When a desingularization of a variety is known, one can study its geometry by 
studying the singular fibres. In our case, we introduce a stratification in 
section 3, which can be viewed as an analogue of Spaltenstein's stratification 
(\cite{Spa}). Its geometric and combinatorial aspects are studied in section 
4, where a combinatorial approach to quiver representations and crystal basis, 
developed in \cite{ReCGS}, plays a key role. As a consequence, certain 
important differences to the nilpotent cone case are established.\\[1ex]
Via the Hall algebra approach, the meaning of our geometric constructions for 
quantum groups is discussed in section 5. The main result is a geometric 
interpretation of the monomial bases introduced in \cite{ReFM}.\\[1ex]
In section 6, we develop a first instance of Springer theory in the quiver 
setup: V. Ginzburg's theory of representations of convolution algebras 
(\cite{CG}) is applied to a quiver analogue of the Steinberg triple variety. 
As a result, we get a 'Kazhdan-Lusztig type' statement, relating certain base 
change coefficients in quantum groups to decomposition numbers for 
representations of certain convolution algebras. The nature of these algebras 
is, however, unknown at the moment, and will be explored in future work.\\[3ex]
{\bf Acknowledgements: } Part of this work was done while the author enjoyed a 
stay at the Weizmann Institute of Science. I would like to thank A. Joseph for 
his kind hospitality. I would like to thank K. Bongartz, W. Borho, V. 
Ginzburg, M. H\"arterich and O. Schiffmann for interesting discussions and 
remarks. This paper was written while the author participated in the TMR 
program "Algebraic Lie Representations" (ERB FMRX-CT97-0100).

\section{Construction of desingularizations}\label{const}

The basic setup for the construction of desingularizations is a study of the 
relation between varieties of quiver representations and of flags in graded 
vector spaces. This setup is also central in Lusztig's construction of quantum 
groups using perverse sheaves (see \cite{LuB}); we will review it here, 
thereby fixing the notations. For general facts concerning representations of 
quivers and their geometry, the reader is referred to \cite{Bo}, 
\cite{RiB}.\\[1ex]
Let $Q$ be a finite Dynkin quiver, i.e. an oriented graph with finite sets of 
vertices $I$ and 
arrows $Q_1$, whose underlying unoriented graph is a disjoint union of Dynkin 
diagrams of type $A,D,E$. Let $k$ be a field. For a finite dimensional 
$I$-graded $k$-vector space $V=\oplus_{i\in I}V_i$, we call the formal sum 
$d=\dimv V=\sum_i(\dim_kV_i)i\in{\bf N}I$ the dimension 
vector of $V$; thus, $V\simeq k^d=\oplus_ik^{d_i}$, where $d_i$ denotes the 
$i$-th component of $d$. The vector space $V$ is 
called pure of weight $i$ if $V_j=0$ for all 
$j\not=i$.\\
A pair $({\bf i},{\bf a})$ of finite sequences ${\bf i}=(i_1,\ldots,i_\nu)\in 
I^\nu$, ${\bf 
a}=(a_1,\ldots,a_\nu)\in{\bf N}^\nu$ is called a monomial in $I$; its weight 
is defined  as $d=|({\bf i},{\bf a})|=\sum_ka_ki_k\in{\bf 
N}I$. A flag $F^*$ of $I$-graded subspaces
$$k^d=F^0\supset F^1\supset\ldots\supset 
F^\nu=0$$
is called of type $({\bf i},{\bf a})$ if $F^{k-1}/F^k$ is pure of weight 
$i_k$ and dimension $a_k$, for all $k=1\ldots \nu$.
The set ${\cal F}_{{\bf i},{\bf a}}$ of all flags of type $({\bf
i},{\bf a})$ is acted upon transitively by the group
$G_d=\prod_{i\in I}{\rm GL}(k^{d_i})\subset{\rm GL}(k^d)$. We fix once and for 
all an
arbitrary flag $F_0^*\subset{\cal F}_{{\bf i},{\bf a}}$ and
denote by $P_{{\bf i},{\bf a}}$ its stabilizer under the
$G_d$-action, a parabolic
subgroup of dimension $\sum_{k\leq l\, :i_k=i_l}a_ka_l$. Thus, ${\cal F}_{{\bf 
i},{\bf a}}\simeq G_d/P_{{\bf
i},{\bf a}}$ is a projective algebraic variety of dimension
$\sum_{k<l\, :i_k=i_l}a_ka_l$.\\[1ex]
We denote by $R_d$ the variety
$$R_d=\bigoplus_{\alpha:i\rightarrow j}{\rm
Hom}_k(k^{d_i},k^{d_j})\subset{\rm End}_k(k^d).$$
This affine variety can be viewed on the one hand as the parameter
space for $k$-representations of the quiver $Q$ of dimension
vector $d$, on the other hand as the variety of '$Q_1$-graded'
nilpotent endomorphisms of $k^d$ (the nilpotency follows since $Q$
has no oriented cycles). The action of ${\rm GL}_k(k^d)$ on ${\rm
End}_k(k^d)$ by conjugation restricts to an action of $G_d$ on
$R_d$ via
$$(g_i)_i(M_\alpha)_\alpha(g_i^{-1})_i=(g_jM_\alpha
g_i^{-1})_{\alpha:i\rightarrow j}.$$
The orbits ${\cal O}_M$ for this action are in one-to-one
correspondence with the isomorphism classes $[M]$ of quiver
representations of dimension vector $d$. Denote by  $\bmod kQ$ the category of 
finite dimensional $k$-representations of $Q$. By Gabriel's theorem, the 
isomorphism classes of indecomposable objects $X_\alpha$ in $\bmod kQ$ are in 
bijection with the set $R^+$ of positive roots for the Dynkin diagram 
corresponding to $Q$ via $\dimv X_\alpha=\alpha$, where $R^+$ is identified 
with a subset of ${\bf N}I$ by identifying the simple root to the vertex $i\in 
I$ with $i\in{\bf N}I$. Thus, the set $[\bmod kQ]$ of isoclasses in $\bmod kQ$ 
is in bijection with the set of functions $R^+\rightarrow{\bf N}$. 
Consequently, the number of $G_d$-orbits in $R_d$ is finite for all $d\in{\bf 
N}I$.\\[1ex]
Define $X_{{\bf
i},{\bf a}}$ as the set of pairs $(M,F^*)\in R_d\times{\cal
F}_{{\bf i},{\bf a}}$ such that $M$ and $F^*$ are compatible,
i.e. $M(F^k)\subset F^k$ for all $k=0\ldots\nu$, where $M$ is
viewed as an endomorphism of $k^d$.
Let $Y_{{\bf i},{\bf a}}$ be the subset of $R_d$ consisting of
those $M$ which are compatible with the fixed flag $F_0^*$. It
is easy to see that $Y_{{\bf i},{\bf a}}$ is a $k$-subspace of
$R_d$ of dimension $\sum_{k<l\, :i_k\rightarrow i_l}a_ka_l$.
Using this notation, we can identify $X_{{\bf i},{\bf a}}$ with
the associated fibre bundle $G_d\times^{P_{{\bf i},{\bf
a}}}Y_{{\bf i},{\bf a}}$ via the map
$$\overline{(g,M)}\mapsto (gMg^{-1},gF_0^*).$$
From this observation we get immediately the following:
\begin{enumerate}
\item $X_{{\bf i},{\bf a}}$ is an irreducible smooth algebraic variety 
of
dimension
$$\dim X_{{\bf i},{\bf
a}}=\sum_{{k<l}\atop{i_k=i_l}}a_ka_l+\sum_{{k<l}\atop{i_k\rightarrow
i_l}}a_ka_l,$$
\item the map $X_{{\bf i},{\bf a}}\rightarrow {\cal F}_{{\bf
i},{\bf a}}$ is a homogeneous vector bundle with typical fibre
$Y_{{\bf i},{\bf a}}$,
\item the canonical map $\pi_{{\bf i},{\bf a}}:X_{{\bf i},{\bf
a}}\rightarrow R_d$ is a projective morphism (a collapsing of a
homogeneous bundle in the sense of G. Kempf).
\end{enumerate}
By the definitions, the image of $\pi_{{\bf i},{\bf a}}$ consists
of those representations $M$ of $Q$ which possess a chain of
subrepresentations
$$M=M^0\supset M^1\supset\ldots\supset M^\nu=0$$
such that for all $k=1\ldots\nu$, the subquotient $M^{k-1}/M^k$ is
isomorphic to $E_{i_k}^{a_k}$, the $a_k$-fold direct sum of the
simple object $E_{i_k}\in\bmod kQ$ (since this is the only
representation of dimension vector $a_ki_k$).\\[1ex]
Since $\pi_{{\bf i},{\bf a}}$ is projective, its image is closed.
It is irreducible since $X_{{\bf i},{\bf a}}$ is, and obviously
$G_d$-stable; thus $\pi_{{\bf i},{\bf a}}(X_{{\bf i},{\bf a}})$ is
one orbit closure $\overline{{\cal O}_M}$. The fibre of $\pi_{{\bf
i},{\bf a}}$ over $M$ identifies with the set of flags of type $({\bf 
i},{\bf
a})$ which are compatible with $M$. These facts suggest the following strategy 
for producing a desingularization of $\overline{{\cal
O}_M}$: we need to find a monomial $({\bf i},{\bf a})$ such that $M$
possesses a unique filtration (as above) of type $({\bf i},{\bf
a})$, and is the generic representation with this property.\\[1ex]
To do this, we translate the construction of monomial bases for
quantized enveloping algebras from \cite{ReFM} into the present
geometric setup; the relation between the algebraic approach of \cite{ReFM} 
and the present geometric approach is explained in section \ref{rqe}. For the 
reader's convenience,
we recall the main steps of this construction.
We start by recalling the concept of a directed partition of $R^+$.
\begin{definition} A partition ${\cal I}_*=({\cal I}_1,\ldots,{\cal 
I}_s)$, where $R^+={\cal
I}_1\cup\ldots\cup{\cal I}_s$ is called directed if
\begin{enumerate}
\item ${\rm Ext}^1(X_\alpha,X_\beta)=0$ for all
$\alpha,\beta\in{\cal I}_t$ for $t=1\ldots s$,
\item ${\rm Hom}(X_\beta,X_\alpha)=0={\rm
Ext}^1(X_\alpha,X_\beta)$ for all $\alpha\in{\cal I}_t$,
$\beta\in{\cal I}_u$, $t<u$.
\end{enumerate}
\end{definition}

\remarks $\;$
\begin{enumerate}
\item Using the fact that the category $\bmod kQ$ is
representation-directed (see \cite{RiB}), we see that directed partitions 
do always
exist.
\item The above conditions can also be written in purely
root-theoretic terms using the non-symmetric bilinear form on ${\bf N}I$ 
defined by
$$\langle i,j\rangle=\left\{\begin{array}{ccc}1&,&i=j,\\
-1&,&i\rightarrow j\mbox{ in $Q$},\\
0&,&\mbox{ otherwise,}\end{array}\right.$$
which fulfills
$$\langle\dimv M,\dimv N\rangle=\dim{\rm Hom}(M,N)-\dim{\rm Ext}^1(M,N).$$
Namely, the definition is equivalent to:
\begin{enumerate}
\item $\langle \alpha,\beta\rangle\geq 0$ for all
$\alpha,\beta\in{\cal I}_t$, $t=1\ldots s$,
\item $\langle \alpha,\beta\rangle\geq 0\geq \langle
\beta,\alpha\rangle$ for all $\alpha\in{\cal I}_t$, $\beta\in{\cal
I}_u$, $t<u$.
\end{enumerate}
This can be derived from the following facts (see \cite{RiB}): there
exists a partial ordering $\preceq$ on the indecomposable
representations of $Q$ such that the non-vanishing of ${\rm
Hom}(U,V)$ already implies $U\preceq V$, and we have $\dim{\rm 
Ext}^1(U,V)=\dim {\rm Hom}(V,\tau
U)$ by the Auslander-Reiten formula, where $\tau$ denotes the Auslander-Reiten 
translation on $\bmod kQ$.
\end{enumerate}

We will now construct a pair consisting of a sequence ${\bf i}$ and an additive
function ${\bf a}:[\bmod kQ]\rightarrow{\bf N}^\nu$ from a directed partition 
${\cal I}_*=({\cal I}_1,\ldots,{\cal
I}_s)$, which will be fixed from now on. This pair will be called the monomial 
function associated to ${\cal I}_*$. We choose a total ordering on $I$ such 
that the existence
of an arrow $i\rightarrow j$ implies $i<j$ (which is possible
since there are no oriented cycles in $Q$). For $t=1\ldots s$, we
define a sequence $\omega_t$ in $I$ by writing the elements of the
set
$$\{i\in I\, :\, \alpha_i\not=0\mbox{ for some }\alpha\in{\cal
I}_t\}$$
in ascending order with respect to the chosen ordering on $I$. We
define the sequence ${\bf i}$ as the concatenation
${\bf i}=\omega_1\ldots\omega_s$. To construct the function ${\bf a}:[\bmod 
kQ]\rightarrow{\bf N}^\nu$ on a 
given representation $M$, we define $M_{(t)}$ as the direct sum of all direct 
summands $U$ of $M$ which are isomorphic to some $X_\alpha$ for 
$\alpha\in{\cal I}_t$, for $t=1\ldots s$. Thus, we arrive at a decomposition 
$$M=M_{(1)}\oplus\ldots\oplus M_{(s)}.$$
For each $t=1\ldots s$, we write the word $\omega_t$ as $(i_1\ldots i_u)$ and 
define the sequence ${\bf a}_t(M)$ as
$$(\dimv_{i_1}M_{(t)},\ldots,\dim_{i_u}M_{(t)}).$$
The sequence ${\bf a}(M)$ is given by concatenation: ${\bf a}(M)=({\bf 
a}_1(M),\ldots,{\bf a}_s(M))$. It is immediately clear that the function ${\bf 
a}$ is 
additive, i.e. ${\bf a}(M\oplus N)={\bf a}(M)+{\bf a}(N)$
(componentwise addition). We can now formulate the main theorem of this 
section.
\begin{theorem} For each representation $M$, 
the map
$$\pi_M=\pi_{{\bf i},{\bf a}(M)}:X_M=X_{{\bf i},{\bf a}(M)}\rightarrow R_d$$
has the orbit closure $\overline{{\cal O}_M}$ as its image, and it restricts 
to an isomorphism $\pi_M:\pi_M^{-1}({\cal 
O}_M)\simeq{\cal O}_M$ over the orbit ${\cal O}_M$. 
Thus, $\pi_M$ is a desingularization of the orbit closure 
$\overline{{\cal O}_M}$.
\end{theorem}

\proof For the calculation of the image of $\pi_M$, we use 
Lemma 4.4 of \cite{ReFM} (see also Proposition 2.4 of \cite{ReGE}), which 
states 
in 
particular the following:\\[1ex]
Let $M$ and $N$ be representations of $Q$ such that ${\rm Ext}^1(M,N)=0$. Let 
$M'$ and $N'$ be degenerations of $M$ and $N$, respectively, and let $X$ be 
an extension of $M'$ by $N'$. Then $X$ is a degeneration of $M\oplus N$. If, 
moreover, ${\rm Hom}(N,M)=0$ and $X\simeq M\oplus N$, then $M'\simeq M$ and 
$N'\simeq N$.\\[1ex]
Here, a representation $M$ is said to degenerate to $N$ if the 
point $N\in R_d$ belongs to the orbit closure $\overline{{\cal O}_M}$; in this 
case we write $M\leq N$. By repeated application, this statement immediately 
generalizes to the 
following one:
\begin{lemma}\label{(*)} If $M_{(1)},\ldots,M_{(s)}\in\bmod kQ$ satisfy ${\rm 
Ext}^1(M_{(k)},M_{(l)})=0$ for all $k<l$, and if $X$ is a representation 
possessing a filtration
$$X=X^0\supset X^1\supset\ldots\supset X^s=0$$
such that $X^{t-1}/X^t$ is a degeneration of $M_{(t)}$ for $t=1\ldots s$, then 
$X$ is a degeneration of the direct sum $M_{(1)}\oplus\ldots\oplus 
M_{(s)}$. If, moreover, ${\rm Hom}(M_{(k)},M_{(l)})=0$ for $k>l$ and $X\simeq 
M_{(1)}\oplus\ldots\oplus M_{(s)}$, then $X^{t-1}/X^t\simeq M_{(t)}$ for all 
$t=1\ldots s$.\end{lemma}
Now suppose a representation $X$ belongs to the image of $\pi_M$. By 
definition, $X$ has a filtration
$$X=X_0\supset X_1\supset\ldots\supset X_s=0$$
such that $$X_{k-1}/X_k\simeq E_{i_k}^{{\bf a}(M)_k}\mbox{ for }k=1\ldots\nu.$$
By writing ${\bf i}$ and ${\bf a}(M)$ in blocks ${\bf 
i}=\omega_1\ldots\omega_s$, ${\bf a}(M)=({\bf a}_1(M),\ldots,{\bf a}_s(M))$ as 
in their above definition, we can 
coarsen this filtration to another one
$$X=X^0\supset X^1\supset\ldots\supset X^s=0$$
such that $X^{t-1}/X^t$ is of dimension vector 
$$\sum_{v=1}^u\dimv_{i_v}M_{(t)}=\dimv M_{(t)}\mbox{ for }t=1\ldots s,$$
where 
$\omega_t=(i_1\ldots i_u)$. Since ${\rm Ext}^1(X_\alpha,X_\beta)=0$ for all 
$\alpha,\beta\in{\cal I}_t$ by the definition of a directed partition, we 
have ${\rm Ext}^1(M_{(t)},M_{(t)})=0$, which means that the orbit of $M_{(t)}$ 
is dense in its variety of representations. Thus, each representation of 
dimension 
vector $\dimv M_{(t)}$ is a degeneration of $M_{(t)}$, which applies in 
particular to the $X^{t-1}/X^t$. Using again the definition of a directed 
partition, we have ${\rm Ext}^1(M_{(t)},M_{(u)})=0$ for all $t<u$. Thus, 
we can apply Lemma \ref{(*)} to conclude that $X$ is a degeneration of 
$M_{(1)}\oplus\ldots\oplus M_{(s)}=M$, which proves that the image of 
$\pi_M$ is contained on $\overline{{\cal O}_M}$. Since we 
already know that the image of $\pi_M$ is an orbit closure, 
we only need to show that $M$ itself belongs to it. We have a filtration
$$M=M_{(1)}\oplus\ldots\oplus M_{(s)}\supset M_{(2)}\oplus\ldots\oplus 
M_{(s)}\supset\ldots\supset M_{(s)}\supset 0$$
by definition. On the other hand, it is easy to see that each representation 
$X$ of $Q$ has a unique filtration $$X=X_0\supset X_1\supset\ldots\supset 
X_m=0$$ such that $$X_{p-1}/X_p\simeq E_{i_p}^{b_p}\mbox{ for }p=1\ldots m,$$
where $i_1\ldots i_m$ is an enumeration of the vertices $i_p\in I$ such that 
$b_p=\dimv_{i_p}X\not=0$, ascending with respect to the chosen ordering on 
$I$. 
In particular, this can be applied to the subfactors $M_{(k)}$ of the 
filtration above. Refining this filtration in this way, it is clear from the 
definitions that we arrive at a filtration of $M$ of type $({\bf i},{\bf 
a}(M))$, proving the first part of the theorem.\\[1ex]
To prove the second part, it suffices to show that $\pi_M$ 
is bijective over ${\cal O}_M$: indeed, the orbit ${\cal O}_M$ is isomorphic 
to the homogeneous space $G_d/{\rm Aut}_{kQ}(M)$, and the restriction of 
$\pi_M$ 
to 
the inverse image of the orbit is thus the projection from an associated fibre 
bundle to a homogeneous space. So assume we are 
given a filtration of $M$ of type $({\bf i},{\bf a}(M))$. Coarsening it as 
above, we obtain a filtration
$$M=M^0\supset M^1\supset\ldots\supset M^s=0$$
such that $M^{t-1}/M^t$ is of dimension vector $\dimv M_{(t)}$, i.e. a 
degeneration of $M_{(t)}$, for $t=1\ldots s$. Since ${\rm 
Hom}(M_{u)},M_{(t)})=0$ for $t<u$ by the properties of a directed partition, 
we can apply the second part of Lemma \ref{(*)} to conclude that 
$M^{t-1}/M^t\simeq M_{(t)}$ for $t=1\ldots s$. Using the ${\rm Ext}$-vanishing 
property for the $M_{(t)}$, an easy induction yields $$M^t\simeq 
M_{(t+1)}\oplus\ldots\oplus M_{(s)}\mbox{ for }t=1\ldots s.$$
But by the ${\rm 
Hom}$-vanishing property, this filtration is unique. As already remarked 
above, the intermediate steps of the original filtration of $M$ are also 
unique 
for general reasons. Thus, the whole filtration is unique, proving the second 
part of the theorem. \hb

\examples As a first example, we show that our desingularizations include as a 
special case the desingularization for equioriented quivers of type $A$ from 
\cite{ADK}. So let $Q$ be the quiver $1\rightarrow 
2\rightarrow\ldots\rightarrow n$. It is known (compare \cite{ReFM}, section 
8) that
$${\cal I}_t=\{\alpha_{n+1-t}+\ldots+\alpha_k\, :\, k=n+1-t,\ldots ,n\}$$ 
defines
a directed partition of $R^+=\{\alpha_i+\ldots+\alpha_j\, :\, 1\leq i\leq 
j\leq n\}$, where $\alpha_i$ denotes the simple root corresponding to the 
vertex $i=1\ldots n$. We denote by $m_{ij}$ the multiplicity of the 
indecomposable representation $E_{ij}=X_{\alpha_i+\ldots+\alpha_j}$ in a 
representation $M$ for $1\leq i\leq j\leq n$. Then the partition ${\cal 
I}_*$ induces the monomial function given by
$${\bf i}=(n,\; n-1,n\; ,\ldots,\; i,\ldots, n\; ,\ldots,\; 1,\ldots, n)$$
and
$${\bf 
a}(M)=(m_{nn},\; m_{n-1,n-1}+m_{n-1,n},m_{n-1,n}\; ,\ldots,\; 
m_{ii}+\ldots+m_{in},\ldots,m_{in}\; ,\ldots$$
$$\ldots,\; m_{11}+\ldots+m_{1n},\ldots,m_{1n}).$$
A short calculation shows that ${\cal F}_M={\cal F}_{{\bf i},{\bf 
a}(M)}$ thus consists of tuples of flags
$$(k^{d_1}=F_1^0\supset F_1^1=0,k^{d_2}=F_2^0\supset F_2^1\supset 
F_2^2=0,\ldots$$
$$,\ldots k^{d_n}=F_n^0\supset F_n^1\supset\ldots\supset F_n^n=0)$$
such that $$\dim F_i^t=\sum_{j+t\leq i\leq k}m_{jk}.$$ 
The variety $X_M$ consists of tuples 
$$(C_1,\ldots,C_{n-1},F_*^*)\mbox{ such that }C_iF_i^*\subset 
F_{i+1}^*\mbox{ for all }i=1\ldots n-1,$$
where $C_i$ is a $d_{i+1}\times 
d_i$-matrix for $i=1\ldots n-1$, and $F_*^*$ is a tuple of flags as 
above. Since the above formula for $\dim F_i^t$ can be interpreted as
$$\dim F_i^t={\rm rank}(C_{i-1}\cdot\ldots\cdot C_{i-t}),$$
this is precisely 
the definition of \cite{ADK}.\\[1ex]
Desingularizations for the same quiver from other directed partitions involve 
compatibility conditions between certain kernels and images of products of the 
matrices $C_i$; they give an explicit realization of the 'mixed' 
desingularizations predicted in \cite{ZePA}.\\[1ex]
As a second example, we sketch a construction of a desingularization for the 
orbit closures of the diagonal action of ${\rm GL}_N$ on a triple product of 
Grassmanians
$${\rm GL}_N:X={\rm Gr}_{d_1}^N\times{\rm Gr}_{d_2}^N\times{\rm Gr}_{d_3}^N$$
for integers $0\leq d_1,d_2,d_3\leq N$; the details will be left to the 
reader. We consider the quiver $Q$ of type $D_4$ with set of vertices 
$I=\{0,1,2,3\}$ and arrows $\alpha_i$ pointing from $i$ to $0$ for $i=1,2,3$, 
respectively. We define a dimension vector $d=N\cdot 0+d_1\cdot 1+d_2\cdot 
2+d_3\cdot 3\in{\bf N}I$ and consider the open subvariety $R_d^0$ of the 
representation variety $R_d$ consisting of tuples $(M_{\alpha_i})$ such that 
all $M_{\alpha_i}$ are injective maps. The subgroup ${\rm GL}_{d_1}\times{\rm 
GL}_{d_2}\times{\rm GL}_{d_3}$ of $G_d$ acts freely on $R_d^0$ with geometric 
quotient $X$. The $G_d$-action on $R_d$ induces the diagonal ${\rm 
GL}_N$-action on $X$, so that we can work with the action of $G_d$ on $R_d^0$ 
equally. Using this relation, the orbits can be described as in (\cite{Ka}, 
1.17.). In particular, the multiplicities of the indecomposables $X_\alpha$ in 
a 
representation $M$ corresponding to a triple of subspaces can be described in 
terms of dimensions of sums and intersections of the $U_i$. A directed 
partition of $R^+$ is given by
$${\cal I}_1=\{\alpha_0\},\; {\cal 
I}_2=\{2\alpha_0+\sum_{i=1}^3\alpha_i,\alpha_0+\alpha_i\, :\, i=1,2,3\},$$
$${\cal 
I}_3=\{\alpha_0+\sum_{i=1}^3\alpha_i,\alpha_0+\sum_{i=1}^3\alpha_i-\alpha_j\, 
:\, j=1,2,3\},\; {\cal I}_4=\{\alpha_i\, :\, i=1,2,3\}.$$
Using this information, one can construct the desingularization of 
$\overline{{\cal O}_M}$ and translate it to the product of Grassmanians, 
yielding the following result:\\
A desingularization of the ${\rm GL}(N)$-orbit closure of a triple 
$(U_1,U_2,U_3)$ is 
given by the projection to the first three factors of the smooth variety
$$\{(U_1,\ldots,U_8)\in\prod_{i=1}^8{\rm Gr}_{d_i}^N\, :$$
$$U_{i+3}\subset U_i\mbox{ for }i=1,2,3,\; U_1+U_2+U_3\subset U_7,\; 
U_4+U_5+U_6\subset U_8,\; U_8\subset U_7\},$$
where $\dim U_{i+3}=\sum_{j\not=i}\dim U_i\cap U_j-\dim U_1\cap U_2\cap 
U_3\mbox{ 
for }i=1,2,3$, $\dim U_7=\dim U_1+U_2+U_3,\; \dim U_8=\sum_{j\not=k}\dim 
U_j\cap U_k-\dim 
U_1\cap U_2\cap U_3$.\\[3ex]
Returning to the general case, we briefly discuss the following:\\
Since there 
are many possible directed partitions for a given quiver $Q$, we 
have constructed a whole class of desingularizations of the orbit closures. 
The question arises how these different desingularizations are related.\\[1ex]
Call a directed partition ${\cal I}'_*$ a refinement of another one ${\cal 
I}_*$ if each ${\cal I}_t$ is a union of several ${\cal I}'_u$. Given such 
partitions ${\cal I}_*, {\cal I}'_*$, let $({\bf i},{\bf a}={\bf a}(M))$ 
(resp. 
$({\bf i'},{\bf a'}={\bf a'}(M))$) be the resulting monomial functions. It is 
easy to see that a flag in ${\cal F}_{{\bf i'},{\bf a'}}$ is a refinement of a 
flag in ${\cal F}_{{\bf i},{\bf a}}$. Thus, we have a canonical projection 
${\cal F}_{{\bf i'},{\bf a'}}\rightarrow{\cal F}_{{\bf i},{\bf a}}$, an 
associated inclusion $P_{{\bf i'},{\bf a'}}\rightarrow P_{{\bf i},{\bf a}}$ 
and an inclusion $Y_{{\bf i'},{\bf a}'}\rightarrow Y_{{\bf i},{\bf a}}$. Using 
the realization of $X_ {{\bf i},{\bf a}}$ as $G_d\times^{P_{{\bf i},{\bf 
a}}}Y_{{\bf i},{\bf a}}$, we see that we arrive at the following statement:
\begin{lemma} There exists a projection $X_{{\bf 
i'},{\bf a'}}\rightarrow X_{{\bf i},{\bf a}}$, which is compatible with the 
maps 
$\pi_{{\bf i'},{\bf a'}}$ and $\pi_{{\bf i},{\bf a}}$.
\end{lemma}

\remark The desingularization given by a directed partition ${\cal I}_*$ is 
thus 
'optimal' if ${\cal I}_*$ is as coarse as possible. This leads to the notion 
of a regular partition in the sense of \cite{ReFM}.

\section{Construction of stratifications}

In this section, we construct a decomposition of the fibres of the maps 
$\pi_{{\bf i},{\bf a}}$ and of the so-called orbital varieties $Y_{{\bf 
i},{\bf a}}\cap{\cal O}_N$ into 
disjoint, irreducible, smooth subvarieties. 
We will use the term stratification in this broad sense. In fact, this 
construction works for arbitrary monomials $({\bf i},{\bf a})$, not 
only for the ones constructed in the previous section. This special case will 
be studied in more detail at the end of the following section. So fix a 
monomial $({\bf 
i},{\bf a})$ and denote by ${\cal F},P,Y,X,\pi$ the objects ${\cal F}_{{\bf 
i},{\bf a}},P_{{\bf i},{\bf a}},\ldots$. First we show that it suffices to 
stratify the orbital varieties $Y\cap{\cal O}_N$.

\begin{lemma}\label{l41} For all $N\in R_d$, we have $G_d$-equivariant 
isomorphisms:
$$G_d\times^{{\rm Aut}_{kQ}(N)}\pi^{-1}(N)\simeq\pi^{-1}({\cal O}_N)\simeq 
G_d\times^P(Y\cap{\cal O}_N).$$
\end{lemma}

\proof The restriction $\pi|_{\pi^{-1}{\cal O}_N}:\pi^{-1}{\cal 
O}_N\rightarrow {\cal O}_N$ is a $G_d$-equivariant projection onto the 
homogeneous space ${\cal O}_N\simeq G_d/{\rm Aut}_{kQ}(N)$, yielding the first 
isomorphism. For the second one, we interpret the map $\pi$ as
$$\pi:G_d\times^PY\rightarrow R_d,\;\; \overline{(g,y)}\mapsto gyg^{-1}.$$
Thus, the subvariety $\pi^{-1}{\cal O}_N\subset R_d$ identifies with
$$\{\overline{(g,y)}\, :\, gyg^{-1}\in{\cal O}_N\}=\{\overline{(g,y)}\, :\, 
y\in{\cal O}_N\}=G_d\times^P(Y\cap{\cal O}_N).$$ \hb

Using these isomorphisms, we can translate geometric properties between the 
fibres of $\pi$ and the orbital varieties. In particular, we can define our 
wanted stratification for $Y\cap{\cal O}_N$, which will simplify the proofs. 
\\[1ex]
We define the stratification in terms of the representation theory of $Q$ 
first:
\begin{definition} Choose an arbitrary sequence $(N=[N_0],\ldots,[N_\nu]=0)$ 
of isomorphism 
classes of 
representations of $Q$ and define
$${\cal S}_{[N_*]}={\cal S}_{([N_0],\ldots,[N_\nu])}=\{L\in Y\, :\, 
L|_{F_0^k}\simeq N_k\mbox{ for all }k=0\ldots\nu\}.$$
\end{definition}

This is well-defined since $Y$ consists of those representations which are 
compatible with the flag $F_0^*$; in particular, each $L|_{F_0^k}$ for 
$L\in Y$ is again a representation of $Q$. Obviously, we have 
$$Y=\bigcup_{N_*}{\cal S}_{[N_*]}\mbox{ and }Y\cap{\cal 
O}_N=\bigcup_{[N_*]:N_0\simeq N}{\cal S}_{[N_*]}.$$
It is also clear that
$${\cal S}_{[N_*]}=\emptyset\mbox{ unless }\dimv 
N_k=d^k:=\sum_{l>k}a_li_l\in{\bf N}I\mbox{ for all }k=0\ldots\nu.$$
Since there are only finitely many isoclasses of representations of a fixed 
dimension vector, we see that the above decomposition is finite.

\begin{proposition} The strata ${\cal S}_{[
N_*]}$ are irreducible, smooth subvarieties of $Y\cap{\cal O}_N$.
\end{proposition}

\proof First, we will derive a geometric construction of the strata.
Denote by $Y_1$ the subspace of $R_d$ of representations which are compatible 
with $F_0^1$ only:
$$Y_1=\{N\in R_d\, :\, N(F_0^1)\subset F_0^1\}.$$
The natural projection $p:Y_1\rightarrow R_{d^1}$, which maps $N\in Y_1$ to 
the 
restriction $N|_{F_0^1}$, is a trivial vector bundle. It restricts to a 
projection $\overline{p}:Y_1\cap{\cal O}_N\rightarrow R_{d^1}.$ It is easy to 
see from the definitions that
$${\cal S}_{[N_*]}=\overline{p}^{-1}({\cal S}_{[N_*']}),$$
where $[N_*']$ denotes the truncated sequence $([N_1],\ldots,[N_\nu])$.
Thus, we can construct the stratum ${\cal S}_{[X_*]}$ inductively, starting 
from ${\cal S}_{([X_\nu])}=R_{d^k}=R_0=0$.\\[1ex]
We will now study the morphism $\overline{p}:Y_1\cap{\cal O}_N\rightarrow 
R_{d^1}$ in more detail. The stratum ${\cal S}_{[N_*']}$ is contained in 
${\cal 
O}_{N_1}$, thus we only need to study the induced morphism 
$\overline{p}:\overline{p}^{-1}{\cal O}_{N_1}\rightarrow{\cal O}_{N_1}$. The 
group $G_{d^1}$ acts on ${\cal O}_{N_1}$, and on $\overline{p}^{-1}{\cal 
O}_{N_1}$ as a subgroup of the parabolic in $G_d$ leaving $F_0^1$ stable. 
Since $\overline{p}$ is obviously equivariant for these actions, we have 
$$\overline{p}^{-1}({\cal 
O}_{N_1})\simeq G_{d^1}\times^{{\rm Aut}_{kQ}(N_1)}\overline{p}^{-1}(N_1).$$
Therefore, it suffices to study the fibre $\overline{p}^{-1}(N_1)$, which 
consists of all points $y$ of ${\cal O}_N$ which are compatible with $F_0^1$ 
and restrict to $N_1$ on this subspace. The inverse image of 
$\overline{p}^{-1}(N_1)$ under the quotient map $$q:G_d\rightarrow{\cal 
O}_N,\; 
q(g)=g^{-1}Ng$$
equals the set of all $g\in G_d$ such that $g^{-1}Ng$ is 
compatible with $F_0^1$ and restricts to $N_1$ on this subspace. An easy 
calculation shows that these conditions are equivalent to $Nf=fN_1$, where $f$ 
denotes the induced homomorphism $f=g|_{F_0^1}:F_0^1\rightarrow 
k^d$. On the other hand, we can consider the morphism
$$pr:G_d\rightarrow {\rm IHom}_k(F_0^1,k^d)=\{f\in{\rm Hom}_k(F_0^1,k^d)\, :\, 
f\mbox{ injective}\}$$
given by restriction of an automorphism of $k^d$ to $F_0^1$,
which is a flat morphism since it is open in the trivial vector bundle ${\rm 
End}_k(k^d)\rightarrow{\rm Hom}(F_0^1,k^d)$. The variety ${\rm 
IHom}_k(F_0^1,k^d)$ contains the subvariety of injective homomorphisms of 
$Q$-representations ${\rm IHom}_{kQ}(N_1,N)$, and its inverse image under $pr$ 
is precisely $q^{-1}\overline{p}^{-1}(N_1)$, as calculated above. Thus, the 
varieties ${\rm IHom}_{kQ}(N_1,N)$ and $\overline{p}^{-1}(N_1)$ are related 
(via $q^{-1}\overline{p}^{-1}(N_1)$) by a quotient morphism on one side and a 
flat morphism on the other side. Since ${\rm IHom}_{kQ}(N_1,N)$ is an open 
subset of the vector space ${\rm Hom}_{kQ}(N_1,N)$, it is an irreducible 
smooth variety. We conclude that the same holds for $\overline{p}^{-1}(N_1)$.
By induction, we see that the proposition is proved. \hb

The above construction also allows us to compute the dimension of each stratum 
${\cal S}_{[N_*]}$.

\begin{proposition}\label{sdim} If the stratum ${\cal S}_{[N_*]}$ is 
non-empty, then it is 
of dimension
$$\dim{\cal S}_{[N_*]}=\dim P_{{\bf i},{\bf a}}+\sum_{k=1}^\nu 
(\dim{\rm Hom}(N_k,N_{k-1})-\dim{\rm End}(N_{k-1})).$$
\end{proposition}

\proof The morphism $q$ defined in the previous proof has relative dimension 
$\dim{\rm 
Aut}_{kQ}(N)=\dim{\rm End}(N)$, and the morphism $pr$ has relative dimension
$$\sum_{i\in I}(d-d^1)_id_i=a_1d_{i_1}=\sum_{k=2}^\nu a_1a_k.$$
Thus, we have
\begin{eqnarray*}
\dim \overline{p}^{-1}(N_1)&=&\dim q^{-1}\overline{p}^{-1}(N_1)-\dim{\rm 
End}(N)\\
&=&\underbrace{\dim{\rm IHom}_{kQ}(N_1,N)}_{=\dim{\rm Hom}(N_1,N)}-\dim{\rm 
End}(N)+\sum_{k=2}^\nu a_1a_k.\end{eqnarray*}

By induction, this yields the formula stated in the proposition. \hb

Concerning the structure of the closure of the strata, we prove, at the 
moment, only the following general criterion. A more detailed study will 
follow in the next section. 

\begin{proposition}\label{sclos} If ${\cal S}_{[L_*]}$ belongs to the closure 
$\overline{{\cal S}_{[N_*]}}$, then $L_k$ belongs to the orbit closure 
$\overline{{\cal O}_{N_k}}$, and ${\cal S}_{[L_k],\ldots,[L_\nu]}$ belongs to 
the closure of ${\cal S}_{[N_k],\ldots,[N_\nu]}$ for all $k=0\ldots\nu$.
\end{proposition}

\proof Suppose ${\cal S}_{[L_*]}\subset\overline{{\cal S}_{[N_*]}}$. 
Considering the $G_d$-saturation of both sides, we get immediately that $N_0$ 
degenerates to $L_0$. Moreover, we have a chain of inclusions
$${\cal S}_{[L_*']}=\overline{p}({\cal S}_{[L_*]})\subset 
\overline{p}(\overline{{\cal S}_{[N_*]}})\subset\overline{\overline{p}({\cal 
S}_{[N_*]})}=\overline{{\cal S}_{[N_*']}}.$$
In this chain, the equalities follow from the properties of the morphism 
$\overline{p}$, and the inclusions are obvious. By induction, the proposition 
follows. \hb

\section{Geometry and combinatorics of the stratification}

We continue to use the notations of the previous section. Using the relation 
between the fibres of $\pi$ and the orbital varieties 
established before, we will now translate the results of the previous section 
to study the geometry of the singular fibres.\\[1ex]
Recall from Lemma \ref{l41} the isomorphism
$$G_d\times^{{\rm 
Aut}_{kQ}(N)}\pi^{-1}(N)\simeq 
G_d\times^P(Y\cap{\cal O}_N).$$
By general properties of associated fibre bundles, we can make the following 
definition.

\begin{definition} Given a sequence $[N_*]=([N]=[N_0],\ldots,[N_\nu]=0)$ 
as before, define the stratum ${\cal F}_{[N_*]}$ in $\pi^{-1}(N)$ via the 
above isomorphism by
$$G_d\times^{{\rm Aut}_{kQ}(N)}{\cal F}_{[N_*]}\simeq G_d\times^P{\cal 
S}_{[N_*]}.$$
\end{definition}

The results of the previous section are immediately translated into the 
following:

\begin{theorem}\label{tos} $\;$
\begin{enumerate}
\item The subsets ${\cal F}_{[N_*]}$ for various sequences $[N_*]$ define a 
stratification of $\pi^{-1}(N)$ into irreducible, smooth, locally closed 
subsets.
\item The stratum ${\cal F}_{[N_*]}$ is non-empty if and only if there exist
short exact sequences $\ses{N_k}{N_{k-1}}{E_{i_k}^{a_k}}$ for all 
$k=1\ldots\nu$.
\item If ${\cal F}_{[N_*]}$ is non-empty, then its dimension equals 
$$\dim{\cal F}_{[N_*]}=\sum_{k=1}^\nu(\dim{\rm Hom}(N_k,N_{k-1})-\dim{\rm 
End}(N_k)).$$
\item If the closure of a stratum ${\cal F}_{[L_*]}$ contains another one 
${\cal F}_{[N_*]}$, then for all $k=1\ldots \nu$, we have $L_k\leq N_k$ and
$\overline{{\cal F}_{([L_k],\ldots,[L_\nu])}}\supset{\cal 
F}_{([N_k],\ldots,[N_\nu])}$.
\end{enumerate}
\end{theorem}

From this, we can derive a formula for the dimension of the singular fibres:

\begin{corollary}\label{dimf} For each $N$, we have
$$\dim\pi^{-1}(N)=\max_{[N_*]:{\cal 
F}_{[N_*]}\not=\emptyset}\sum_{k=1}^\nu(\dim{\rm Hom}(N_k,N_{k-1})-\dim{\rm 
End}(N_k)).$$
\end{corollary}

Another application of the stratification concerns the determination of the 
irreducible components of the singular fibres $\pi^{-1}(N)$.

\begin{proposition}\label{irrc} $\;$
\begin{enumerate}
\item If $I$ is an irreducible component of $\pi^{-1}(N)$, then 
there exists a 
sequence $[N_*]$ such that $I=\overline{{\cal F}_{[N_*]}}$.
\item The subvariety $\overline{{\cal F}_{[N_*]}}$ is an irreducible component 
of 
$\pi^{-1}(N)$ provided the following holds: If $[L_*]$ is a sequence such that 
$$L_k\leq N_k\mbox{ and }\dim{\cal 
F}_{([L_k],\ldots,[L_\nu])}\geq\dim{\cal F}_{([N_k],\ldots,[N_\nu])}$$
for all 
$k=1\ldots\nu$, then already $L_k\simeq N_k$ for all $k=1\ldots\nu$.
\end{enumerate}
\end{proposition}

In the case of Springer's resolution, the strata of the singular fibres are 
in bijection to Young tableaux of a fixed shape (compare \cite{Spa}), i.e. to 
paths in the Young 
graph, which has Young diagrams as vertices, and arrows corresponding to 
additions of a single box. Analoguously, we will now construct a graph 
structure on isomorphism 
classes of representations of $Q$ whose paths of a certain type parametrize 
the strata in our stratification.

\begin{definition} Let $\Gamma(Q)$ be the $I\times{\bf N}$-coloured graph with 
vertices the isomorphism classes of representations of $Q$, and with an arrow 
of colour $(i,n)\in I\times{\bf N}$ from $N$ to $X$ if there exists an exact 
sequence $\ses{N}{X}{E_i^n}$. A path of colour $({\bf i},{\bf a})$ in 
$\Gamma(Q)$ is a sequence $[N_*]$ such that there exist arrows 
$[N_{k-1}]\rightarrow[N_k]$ of colour $(i_k,a_k)$ for all $k=1\ldots\nu$.
\end{definition}

The second part of Theorem \ref{tos} now reads as follows:

\begin{lemma}\label{par} The non-empty strata in $\pi_{{\bf i},{\bf 
a}}^{-1}(N)$ are 
parametrized by the 
paths of colour $({\bf i},{\bf a})$ from $[0]$ to $[N]$ in $\Gamma(Q)$.
\end{lemma}

For a certain class of quivers, which we call special, the strata can be 
pa\-ra\-metrized in a completely combinatorial way, using results from 
\cite{ReCGS}. Since this class of quivers was described there only in 
representation-theoretic terms, we first give a purely root-theoretic 
description.

\begin{definition} A quiver $Q$ is called special if for all indecomposables 
$X_\alpha$ and all simples $E_i$, we have
$\dim{\rm Hom}(X_\alpha,E_i)\leq 1$.
\end{definition}

Call a vertex $i\in I$ thick if there exists a root $\alpha\in R^+$ such that 
the simple root $i\in R^+$ appears with multiplicity at least $2$ in $\alpha$.

\begin{proposition} $Q$ is special if and only if no thick vertex is a source 
of $Q$.
\end{proposition}

\proof First, we prove the following statement:\\
$Q$ is special if and only if $\langle\alpha,i\rangle\leq 1$ for 
all 
$\alpha\in R^+$ and all $i\in I$.\\
So suppose $Q$ to be special, and let $\alpha\in R^+$ and $i\in I$. Then
$$\langle\alpha,i\rangle=\dim{\rm Hom}(X_\alpha,E_i)-\dim{\rm 
Ext}^1(X_\alpha,E_i)\leq\dim{\rm Hom}(X_\alpha,
E_i)\leq 1.$$
Conversely, suppose $X_\alpha$ and $E_i$ as above are given. If
${\rm Hom}(X_\alpha,E_i)=0$, there is nothing to prove. Otherwise, we already 
have ${\rm Ext}^1(X_\alpha,E_i)=0$ by the directedness of $\bmod kQ$, so
$$1\geq \langle\alpha,i\rangle=dim{\rm Hom}(X_\alpha,E_i)-\dim{\rm 
Ext}^1(X\alpha,E_i)=\dim{\rm Hom}(X_\alpha,E_i).$$

Now we come to the proof of the proposition. First we show that the above 
condition is 
neccessary. So let $i$ be
a thick vertex of $Q$, and let $\alpha=\sum_{i\in I}\alpha_ii\in R^+$ be a 
root such that $\alpha_i\geq 2$. By the previous 
lemma, we have
$$1\geq \langle\alpha,i\rangle=\alpha_i-\sum_{j\rightarrow i}\alpha_j\geq 
2-\sum_{j\rightarrow i}\alpha_j,$$
so that $i$ cannot be a source in $Q$.\\[1ex]
To prove the converse, we have to proceed by a case-by-case analysis:\\
If $Q$ is of type $A_n$, then $\alpha_i$ equals $0$ or $1$ for all
$\alpha\in R^+$ and $i\in I$. This means that no vertex is thick, and that
$Q$ is always special. If $Q$ is of type $D_n$, we have the following 
property: If $\alpha_i$=2,
then $\alpha_j=1$ for all $j$ connected to $i$. Thus, we are done by the
above inequality, assuming that no thick $i$ is a source. If $Q$ is of type 
$E_6$ or $E_7$, then a direct inspection of the set
of positive roots shows the converse. If $Q$ is of type $E_8$, then all 
vertices are thick (look at the
highest root $\alpha=(2,3,4,6,5,4,3,2)$), so there is no orientation $Q$
for which no thick vertex is a source, and thus there is no special 
orientation. \hb

Assume now that $Q$ is special. For each vertex $i\in I$, define a partially 
ordered set ${\cal P}_i\subset R^+$ as the set of all positive roots $\alpha$ 
such that $\langle\alpha,i\rangle=1$, partially ordered by $\alpha\leq\beta$ 
if $\langle \alpha,\beta\rangle\geq 1$ (compare \cite{ReCGS}, Proposition 
4.3). Let ${\cal S}_i$ be the set of antichains in ${\cal P}_i$; it is 
naturally in bijection to the poset of order ideals in ${\cal P}_i$, from 
which 
it inherits a partial ordering. Let $\tau\in W$ be the Weyl group 
element defined by $X_{\tau\alpha}=\tau X_\alpha$. Given an antichain 
$A\in{\cal S}_i$, we 
denote by $l(A)$ the set of all $\alpha\in {\cal P}_i$ which are minimal with 
the property $\{\alpha\}\not\leq A$. Translating Proposition 4.5 of 
\cite{ReCGS} into a root-theoretic language, we get:

\begin{proposition} For any representation $M\in\bmod kQ$, the possible middle 
terms $X$ of exact sequences $\ses{M}{X}{E_i}$ are parametrized by antichains
$A\in{\cal S}_i$ such that $X_{\tau\alpha}$ is a direct summand of $M$ for 
each $\alpha\in l(A)$. The corresponding middle term is given by
$$X=B\oplus\bigoplus_{\alpha\in A}X_\alpha,\mbox{ where 
}M=B\oplus\bigoplus_{\alpha\in l(A)}X_{\tau\alpha}.$$
\end{proposition}

Repeated application of this proposition yields all possible middle terms of 
exact sequences $\ses{M}{X}{E_i^n}$. Combining this with Lemma \ref{par}, we 
get:

\begin{corollary} If $Q$ is special, there exists a purely root-theoretic 
parametrization of the strata in each $\pi_{{\bf i},{\bf a}}(N)$.
\end{corollary}

\remark Using the proof 
of Proposition 5.2 of \cite{ReCGS}, one could also give a combinatorial 
formula for the dimension of the singular fibres.\\[2ex]
\example We consider again the desingularization from the first example of 
section 2, and use the notations defined there. Let $N$ be a point in 
$\overline{{\cal O}_M}$, and assume $N$ is given by multiplicities 
$(n_{ij})_{i\leq j}$. By the above proposition, all exact sequences with 
simple right end term are of the form $0\rightarrow B\oplus 
E_{i+1,j}\rightarrow B\oplus E_{ij}\rightarrow E_i$ for some $j\geq i$. Thus, 
exact sequences with right end term $E_i^n$ are of the form
$$\ses{B\oplus\bigoplus_{j\geq i}E_{i+1,j}^{p_j}}{B\oplus\bigoplus_{j\geq 
i}E_{ij}^{p_j}}{E_i^n}.$$
It follows that the sequences $[N_*]$ defining non-empty strata ${\cal 
F}_{[N_*]}$ in $\pi_M^{-1}(N)$ are parametrized by tuples $(p_k^{ij})$ for 
$i\leq j\leq k$ such that 
$$\sum_kp_k^{ij}=m_{ij}+\ldots+m_{in},\;\; 
\sum_k(p_j^{ik}-p_j^{k,i-1})=n_{ij}\mbox{ for all }i\leq j,$$
$$\sum_{l\leq i}(p_k^{lj}-p_k^{l,j+1})\leq 0\mbox{ for all }i\leq j<k.$$

Now we turn to the special case where $({\bf i},{\bf a})$ is a monomial 
function constructed from some directed partition. In this case, the 
stratification of the singular fibres of $\pi_M$ can be viewed as a graded 
analogue of Spaltenstein's stratification of \cite{Spa}. However, in our 
situation, the geometry of the singular fibres is considerably more 
complicated.\\[2ex]
\remarks $\;$
\begin{enumerate}
\item The fibres are far from being equidimensional, as the dimension formula 
from Theorem \ref{tos} already shows in examples for type $A_3$.
\item The irreducible components are not in bijection to the strata. In fact, 
one frequently encounters inclusions of closures of strata. At the moment, the 
best result for identifying irreducible components is Proposition \ref{irrc}; 
a representation-theoretic criterion for inclusions of closures of strata is 
missing. Nevertheless, in examples of type $A_3$ Proposition \ref{irrc} 
detects roughly half the number of irreducible components.
\item The maps $\pi_M$ are not semi-small (in the sense of intersection 
homology theory) in general, as can 
be seen in 
examples using Corollary \ref{dimf}. In fact, the results of the following 
section relate this property to monomiality of certain elements in Lusztig's 
canonical basis, which is discussed in detail in \cite{LuTM}, \cite{ReMCB}.
\item It is not clear whether the singular fibres are always connected. This 
follows immediately from Zariski's Main Theorem provided the orbit closures 
$\overline{{\cal O}_M}$ are normal varieties. So far, this is only known for 
quivers of type $A$ by \cite{Zwa}, where it is reduced to the case of an
equioriented quiver of type $A$, which was treated in \cite{ADK}. 
\end{enumerate}

In fact, the connectedness of the singular fibres is directly related to the 
question of normality (or at least unibranchness, which, together with a 
Frobenius splitting, implies normality) of the orbit closures, as the 
following lemma shows:

\begin{lemma} Let $f:X\rightarrow Y$ be a desingularization map between 
algebraic 
varieties $X$, $Y$. Then all $f^{-1}(y)$ are connected if and only if $Y$ is 
unibranch, i.e. if the normalization morphism
$g:\widetilde{Y}\rightarrow Y$ is bijective.
\end{lemma}

We omit the proof. It follows easily from the universal property of the 
normalization, standard properties of finite morphism and Zariski's Main 
Theorem.

\section{Relation to quantized enveloping algebras}\label{rqe}

In this section, we assume $k$ to be a finite field with $v^2$ elements.
Given a finite $G$-set $X$, i.e. a finite set on which a group $G$ acts, we 
denote by ${\bf Q}_G[X]$ the set of $G$-invariant functions from $X$ to the 
rationals. We define the Hall algebra of the quiver $Q$ (\cite{RiH}) in terms 
of a convolution product. Denote by $H_v(Q)$ the direct sum
$$H_v(Q)=\bigoplus_{d\in{\bf N}I}{\bf Q}_{G_d}[R_d],$$
where each $R_d$ is viewed as a finite $G_d$-set over the finite field $k$.
Define a multiplication on $H_v(Q)$ by
$$(f*g)(X)=v^{\langle d,e\rangle}\sum_{U\subset X}f(X/U)\cdot G(U),$$
where $f\in{\bf Q}_{G_d}[R_d]$, $g\in{\bf Q}_{G_e}[R_e]$, and the sum runs 
over all subrepresentations of $X$, which is viewed as an object of $\bmod 
kQ$. This 
product is 
well defined since the functions $f$ and $g$ only depend on the isomorphism 
class of a representation; it endows $H_v(Q)$ with the structure of an 
associative, ${\bf N}I$-graded ${\bf Q}$-algebra. Denote by $E_i$ the function 
with value $1\in{\bf 
Q}$ on the one-point set $R_i$. The main result of \cite{RiH} can be stated in 
this language as:\\
The map $\eta:E_i\mapsto E_i$ extends to an isomorphism of 
${\bf N}I$-graded ${\bf Q}$-algebras
$$\eta: H_v(Q)\stackrel{\sim}{\rightarrow}{\cal U}_v({\frak n}^+),$$
where ${\cal U}_v({\frak n}^+)$ denotes the positive part of the quantized 
enveloping algebra to the Dynkin diagram underlying $Q$, which is a ${\bf 
Q}$-algebra with generators $E_i$ for $i\in I$, subjected to the quantized 
Serre-relations at $v$.\\[1ex]
Analoguosly to the previous sections, we fix a directed partition $({\cal 
I}_*)$ of 
$R^+$ and denote by $({\bf i},{\bf a})$ the corresponding monomial function.
Given a representation $M$ of dimension vector $d$, we thus have a monomial 
$({\bf i},{\bf a}=(a_1\ldots a_\nu))$ and the corresponding
desingularization $\pi_M:X_M\rightarrow\overline{{\cal O}_M}$.

\begin{definition} Define a monomial in $H_v(Q)$ by
$$E^{(M)}=E_{i_1}^{(a_1)}*E_{i_2}^{(a_2)}*\ldots 
*E_{i_\nu}^{(a_\nu)}\in{\bf Q}_{G_d}[R_d],$$
where $E_{i}^{(n)}$ denotes the divided power $([n]!)^{-1}E_i^{*n}\in{\bf 
Q}_{G_{ni}}[R_{ni}]$.
\end{definition}

\begin{lemma}\label{mon} For all representation $M$, $N$, we have
$$E^{(M)}(N)=v^{\dim{\rm End}(M,M)-\dim M}|\pi_M^{-1}(N)|,$$
where $|\cdot|$ denotes the cardinality of a finite set.
\end{lemma}

\proof The value of the function $E_i^{(n)}$ on the unique point $0$ of 
$R_{ni}$ is easily computed from the definitions as
$$E_i^{(n)}(0)=([n]!]^{-1}E_i^{*n}(0)=([n]!)^{-1}v^{n(n-1)/2}|{\cal 
F}_n(k)|=v^{n(n-1)},$$
where ${\cal F}_n$ denotes the set of complete flags in the vector space 
$k^n$. Using this, we can compute the value of a monomial in the $E_i^{(n)}$ 
as $E_{i_1}^{(a_1)}*E_{i_2}^{(a_2)}*\ldots 
*E_{i_\nu}^{(a_\nu)}(N)=$
\begin{eqnarray*}
&=&v^{\sum_{k<l}\langle a_ki_k,a_li_l\rangle}\sum_{N=N_0\supset\ldots\supset 
N_\nu=0}E_{i_1}^{(a_1)}(N_0/N_1)\cdot\ldots\cdot 
E_{i_\nu}^{(a_\nu)}(N_{\nu-1}/N_\nu)\\
&=&v^C |\{N_*=(N=N_0\supset\ldots\supset 
N_\nu=0)\, :\, \dimv N_{k-1}/N_k=a_ki_k\mbox{ for all k}\}|\\
&=&v^C|\pi_M^{-1}(N)|
\end{eqnarray*}
where $C=\sum_{k\leq l}\langle a_ki_k,a_li_l\rangle-\sum_ka_k$. Thus, it 
remains to identify the exponents of $v$. We have
$$\dim{\rm End}(M,M)=\dim G_d-\dim{\cal O}_M=\dim G_d-\dim X_M;$$
using the dimension formulae of section \ref{const} and the definition of the 
bilinear form $\langle,\rangle$, this reduces to an easy calculation. \hb

The set of elements $B=\{E_{M}\, :\, [M]\in[\bmod kQ]\}$ of $H_v(Q)$ defined by
$$E_{M}(N)=\left\{\begin{array}{lcc}
v^{\dim{\rm End}(M,M)-\dim M}&,&M\simeq N,\\
0&,&\mbox{ otherwise}
\end{array}\right.$$
is obviously a basis for $H_v(Q)$. It is proved in \cite{L1} that the 
corresponding basis for ${\cal U}_v({\frak n}^+)$ is of PBW type. The above 
lemma can now be rewritten as
\begin{proposition} For each representation $M$, we have
$$E^{(M)}=\sum_{[N]}|\pi^{-1}(N)|E_{N}.$$
In particular, the set of elements ${\cal M}=\{E^{(M)}\, :\, [M]\in[\bmod 
kQ]\}$ is a 
monomial basis for $H_v(Q)$, which has upper unitriangular base change (with 
respect to the degeneration ordering) to the PBW basis $B$.
\end{proposition}

\proof The base change coefficients are precisely calculated by the preceding 
lemma. We have $|\pi^{-1}_M(M)|=1$, and $|\pi^{-1}_M(N)|\not=0$ only if $M\leq 
N$ 
since $\pi_M$ is a desingularization of $\overline{{\cal O}_M}$. This implies 
the unitriangularity property; in particular, the 
$E^{(M)}$ already form a basis for $H_v(Q)$. \hb

\remark It can easily be seen from the definitions that the monomial bases 
thus constructed are exactly the same as in \cite{ReFM}.\\[2ex]
Now we consider the relation of the monomial bases to Lusztig's canonical 
basis. For each isoclass $[M]$ in $\bmod kQ$, we introduce an element ${\cal 
E}_M$ in $H_v(Q)$ by
$${\cal E}_M(N)=v^{\dim{\rm End}(M)-\dim M}\sum_{i\in{\bf Z}}\dim {\cal 
H}^i_N{\cal I\cal C}({\cal 
O}_M)v^i,$$
where ${\cal H}^i_N{\cal I\cal C}({\cal O}_M)$ denotes the stalk at the point 
$N$ of the $i$-th cohomology sheaf of the intersection cohomology complex 
corresponding to the orbit of $M$ as a variety over $k={\bf C}$. By general 
properties of perverse sheaves, 
these elements have upper unitriangular base change (with respect to the 
degeneration ordering) to the PBW basis, thus
${\cal B}=\{{\cal E}_M\, :\, [M]\in[\bmod kQ]\}$
is a basis for $H_v(Q)$.

\begin{proposition} The base change from ${\cal M}$ to ${\cal B}$ is upper 
unitriangular with entries consisting of Laurent polynomials with non-negative 
integer coefficients.
\end{proposition}

\proof The proof uses Lusztig's realization of ${\cal U}_v({\frak n}^+)\simeq 
H_v(Q)$ in terms of convolution of perverse sheaves (\cite{LuB}, 13.2.11.). 
Combining Lemma \ref{mon} with (\cite{LuB}, 9.1.3.) and (\cite{LuB}, 
13.1.12.(b)), one can interprete the base change coefficients from ${\cal M}$ 
to ${\cal B}$ as multiplicities of simple perverse sheaves in the 
decomposition of the derived direct image $(\pi_M)_*1$ of the constant sheaf 
on $X_M$ (up to some shifts). From this, both the positivity property and the 
unitriangularity property follow immediately (since $\pi_M$ is a 
desingularization). \hb

\section{Convolution algebras}

In this section, we assume $k={\bf C}$. Fix ${\cal I}_*$, ${\bf i}$ and ${\bf 
a}$ as in the previous sections. Given a representation $M$, we have the 
corresponding desingularization $\pi_M:X_M\rightarrow\overline{{\cal O}_M}$. 
We define 
$$Z_M:=X_M\times_{\overline{{\cal O}_M}}X_M.$$
This is an analogue of the triple variety in the context of Springer's 
desingularization: Using the description of $X_M$ of section \ref{const}, it 
is easy to 
see that
$$Z_M\simeq\{(N,F_1^*,F_2^*)\in R_d\times{\cal F}_{{\bf i},{\bf a}}\times{\cal 
F}_{{\bf i},{\bf a}}\, :\, NF_i^*\subset F_i^*\mbox{ for }i=1,2\}.$$
We recall some results on convolution algebras from \cite{CG}. The total 
Borel-Moore homology
$$A_M=H_*^{\rm BM}(Z_M,{\bf C})$$
carries a natural associative algebra structure via convolution. It acts via 
convolution on the total Borel-Moore homology $V_N=H_*^{\rm 
BM}(\pi_M^{-1}(N),{\bf C})$ of the fibres of the desingularization $\pi_M$. 
The simple 
$A_M$-modules are parametrized by a subset of the set of isoclasses 
$[N]\in\bmod kQ$ such that $N$ is a degeneration of $M$. In fact, each 
corresponding simple module $L_N$ can be realized as a quotient of $V_N$ (this 
follows from (\cite{CG}, chapter 8), since in our situation, all orbits ${\cal 
O}_N$ are simply connected (the stabilizer ${\rm Aut}_{kQ}(M)$ is conneced). 
The result of \cite{CG} which we need here mainly is 
Theorem 8.6.23, stating that the decomposition numbers of the representations 
$V_N$ are given by Poincare polynomials in local intersection cohomology:
$$[V_{N_1}:L_{N_2}]=\sum_{i\in{\bf Z}}\dim{\cal H}^i_{N_1}{\cal I\cal C}({\cal 
O}_{N_2}),$$
where $[V_{N_1}:L_{N_2}]$ denotes the Jordan-H\"older multiplicity of the 
simple $A_M$-module $L_{N_2}$ in $V_{N_1}$.\\[1ex]
All objects considered in the previous section have 'generic' versions: In the 
${\bf Q}(v)$-algebra ${\cal U}^+$ given by generators $E_i$ for $i\in I$ 
subjected to the quantized Serre relations, one has elements $E_{M},{\cal 
E}_M,E^{(M)}$ which specialize to the corresponding basis elements for each 
$v$ such that $v^2$ is the cardinality of a finite field (\cite{LuB}). In 
particular, the enveloping algebra ${\cal U}({\frak n}^+)$ is provided with 
both canonical and PBW bases.\\[1ex]
This allows us to prove the following 'Kazhdan-Lusztig type' statement:
\begin{theorem} Writing ${\cal E}_M=\sum_{[N]}\zeta^M_NE_{N}$ in the 
enveloping algebra ${\cal U}({\frak n}^+)$, we have
$$\zeta^M_N=[V_N:L_M].$$
\end{theorem}

\proof Since $\pi_M$ is a desingularization, the fibre over $M$ is a single 
point, thus $V_M=L_M$ is already a simple ($1$-dimensional) $A_M$-module. The 
result of \cite{CG} cited above identifies the decomposition number 
$[V_N:L_M]$ with $\sum_{i\in{\bf Z}}\dim{\cal H}^i_N{\cal I\cal C}({\cal 
O}_M)$. But by the definition of the canonical basis given in the previous 
section, this is just the $v=1$-specialization of the base change coefficient 
to the PBW basis. \hb

The immediate problem of an algebraic description of the algebras $A_M$ 
arises. At the moment, no general descriptions are known. The analysis of 
their structure is planned for future publications. At this point, we discuss 
one particular example.\\[2ex]
\example Let $Q$ be the quiver $i\rightarrow j$ of type $A_2$, and let $M$ be 
the representation $E_j^{n-1}\oplus E_{ij}\oplus E_i^{m-1}$, so that 
$\overline{{\cal O}_M}$ is the set of complex $m\times n$-matrices of rank at 
most $1$. Then $A_M$ has generators $x,y$, subjected to the relations
$$x^n=0,\;yx^ky=\delta_{k,n-m-1}\cdot y\mbox{ for all }k\geq 0,\; 
x^m=\sum_{p+q=n-1}x^pyx^q.$$

\end{document}